\documentstyle[11pt,mssymb]{article}

\font \twlmath = msym10 at 11pt
\newfam \mathfam \textfont \mathfam = \twlmath

\newcommand{\arrow}{\rightarrow}
\newcommand{\g}{{\goth g}}
\newcommand{\la}{{\lambda}}
\newcommand{\Z}{{\Bbb Z}}
\newcommand{\C}{{\Bbb C}}

\newcommand{\FF}{{\cal F}}
\newcommand{\I}{{\cal I}}
\newcommand{\LL}{{{\cal L}}}
\newcommand{\OO}{{\cal O}}

\def\square{\hbox{\vrule\vbox{\hrule\phantom{o}\hrule}\vrule}}

\begin{document}

\setlength{\parindent}{7pt}
\setlength{\parskip}{3pt plus 6pt minus 1pt}

\centerline{\bf  KOSZUL PROPERTY AND FROBENIUS SPLITTING}

\centerline{\bf  OF SCHUBERT VARIETIES}
\vskip 1cm
\centerline{\sc Roman Bezrukavnikov}
\vskip 5mm
\centerline {University of Tel-Aviv}

\vskip 1cm

{\bf  Introduction.} Let  $\g$ be a semisimple Lie algebra
over $\C$ and let  $\la$ be a dominant weight.
By $V_\la $ we denote the irreducible representation of  highest weight $\la$.

Consider an associative commutative graded algebra ${\displaystyle
A :=\bigoplus_{n=0}^{\infty} V_{n \la}}$, with  (nonzero) $\g$-invariant
multiplication. (Such multiplication is unique up to isomorphism.)
In other words, $A$ is the projective coordinate ring of the orbit
of the  highest weight vector in projective space  ${\Bbb P} (V_\la^*)$.

A classical theorem of Kostant (see e.g. [LT], Theorem 1.1) asserts
 that algebra $A$ is {\it quadratic}, i.e.
it is generated by its first graded component, and the ideal of
relations is generated
by elements of degree 2. Among  quadratic algebras there is
a class of algebras enjoying  especially
nice homological (and, hence, deformation) properties.
It is the class of {\it Koszul algebras}.
(We refer to [BGS], [PP], [BF] for an exposition of  the theory of
 Koszul algebras).

Below we prove that  $A$, and, more generally, the coordinate ring of
a Schubert variety in any projective embedding belongs to this class.

A method of investigating the coordinate rings of Schubert
varieties, called the Frobenius splitting method,
 was developed in the sequence of papers  [MR], [R] (and reference there).
 Our observation is that this method is
ideally suited for establishing the Koszul property of a commutative ring.

{\bf 1. Frobenius Splitting.} The following definitions
are taken from [MR].

 Let $X$ be a scheme over a
field of characteristic $p>0$. Let $Fr:\;X\arrow X$ be the
(geometric) Frobenius morphism. For a regular function $f$
on an affine open subscheme we have $Fr^*(f)=f^p$. To the morphism
$Fr$ there corresponds the structure morphism of sheaves
$F:\; {\cal O}\arrow Fr_* ({\cal O})$.

{\sf 1.1. Definition.} {\it  $X$ is called Frobenius split if there
exists a morphism $s:\; Fr_* ({\cal O})\arrow
  {\cal O} $ such that  the composition
$s\circ F$ is identity.  Such
 $s$ is called a Frobenius splitting of $X$.}

{\sf 1.2.} Let $Y$ be a closed subscheme of $X$, and $J_Y$ be its ideal sheaf.

{\sf Definition.} {\it A Frobenius splitting $s$ is called
compatible with $Y$ if the image of  $Fr_*(J_Y)$ under
$s$ lies in $J_Y$. If such an $s$ exists, then $Y$
is called compatibly split in $X$.}

{\sf 1.3. Lemma.} (cf. [R], 1.3)
 {\it  i) If a Frobenius splitting is compatible
with a subscheme $Y$ then it induces a Frobenius splitting
of $Y$.

ii) A Frobenius split scheme is reduced.} \hfill \square

{\sf 1.4.} Let $X$ be a projective variety, and ${\cal L}$
be an ample invertible sheaf on $X$.

{\sf Proposition} {\it i) If $X$ is Frobenius split, then
$H^i(X, {\cal L})=0$ for $i>0$.

ii) If $Y$ is compatibly  split in $X$, then
$H^i(X,{\cal L} \otimes J_Y)=0$ for $i>0$. }

{\sf Proof: } see [MR], Proposition 3. \hfill \square

{\sf 1.5.} We shall apply this proposition in the following
situation. Let $G$ be a Chevalley group, $B^+,\; B^-$
be its opposite Borel subgroups, and let $P$ be a parabolic
 subgroup. From now on $X:=G/P$ will denote a flag manifold.
 By a Schubert (respectively, an opposite Schubert)
variety we mean a closure of a $B^+$ (respectively, a $B^-$)
orbit in $X$, and by a relative Schubert variety a closure
of a $G$-orbit in $X\times X$.

All these objects are considered as schemes over $Spec({\Bbb Z})$.
For any field ${\Bbb F}$ we can get the corresponding objects
over $Spec ({\Bbb F})$ by the base change. We omit the base scheme
 in the notations when it causes no confusion.

{\sf 1.6.} Let us fix $n$ and consider the following
subvarieties of $X^n$:

a) preimage of a Schubert variety under the projection on the
first coordinate

b) preimage of a  relative Schubert variety under the projection on the
$(i,i+1)$ coordinates for all $i=1,..,n-1$.

c) preimage of an opposite  Schubert variety under the projection on the
$n$-th coordinate

d) $X^n,\;\emptyset$

Let $S^n$ denote the smallest  set of subschemes of $X^n$,
which contains all the subschemes a),b),c),d) and is closed under the
 operations of union and (scheme-theoretic) interesection of subschemes.

{\sf Proposition.} {\it For  any   $p>0$
  there exists a Frobenius splitting
of $X^n$ over a  field of characteristic $p$, which is compatible
with all subschemes $Z\in S^n$.}

{\sf Proof } will be given in section 3.

{1.7.}  We consider the schemes over the spectrum of
a field. Let ${\cal L}$ be an ample invertible sheaf
on $X^n$.

{\sf Corollary.} {\it i) Any scheme $Z\in S^n$ is reduced.

ii) For $Z_1 \subset Z_2$ put $J^{Z_2}_{Z_1}$ to be the kernel
of the restriction morphism ${\cal O}_{Z_2}\arrow {\cal O}_{Z_1}$.
Then  for all $Z_1,\;Z_2
\in S^n$ we have: $H^i(X^n,J^{Z_2}_{Z_1} \otimes {\cal L})=0$
for $i>0$.}

{\sf Proof:} (cf. [MR], Theorem 3) For fields of finite characteristic the
statements
 follow from Proposition 1.6. together with  Lemma 1.3.ii), and
Proposition 1.4.ii) respectively.

Note  that for any field $k$
the set  $S^n_{k}$ of subschemes in $X^n/ Spec (k)$ can be obtained
from the corresponding set of subschemes in $X^n/ Spec (\Z) $ by the
base change. Thus all the sheaves $J^{Z_2}_{Z_1}$ over $k$ are
the base change of the sheaves over $\Z$. Also there exists
an ample sheaf  $\LL _{\Z}$ on $X_{\Z}$ which gives $\LL$ by the base change.

It is easy to see that there exists a nonempty open
 $U\hookrightarrow Spec(\Z)$ such that
all the sheaves under consideration are flat over $U$. The statement
ii) for the base field of characteristic 0 then follows from the
corresponding statement over  finite
fields by semicontinuity ([EGA], Theorem 7.7.5, page 199).

Let us prove i).
Assume that the nil-radical ${\cal I}_Z \hookrightarrow \OO _Z$ is non-zero.
{}From the statement i) for finite fields it follows,
 by Nakayama lemma, that the support of $\I _Z$ does not intersect the fiber
of $X^n$ over the prime ideal $(p)$ for any $p>0$. This is a contradiction
since the support
of any coherent sheaf is closed.
\hfill \square

{\bf 2. Koszul Property.} We are now in a position to prove
the Koszul property. Let us descibe a way to deduce  combinatorial
information from the geometrical statements 1.7.

{\sf 2.1.} Suppose that  a scheme $Y$,  a flat sheaf ${\cal L}$ on it
and a set $S$ of subschemes of $Y$ are given. About $S$ we
assume that $Y\in S$, and that $S$ is
closed under the operations of scheme-theoretic intersection and union
of subschemes. Let us denote $V:=H^0(Y,{\cal L})$; and let $V_Z$
be the kernel of the restriction map  $V\arrow H^0({\cal L}\otimes
{\cal O}_Z)$.

{\sf Lemma.} {\it If for ${Y,\;S,\; {\cal L}}$ the statements of Corollary 1.7.
are satisfied, then

a)  the restriction map  $V\arrow H^0({\cal L}\otimes
{\cal O}_Z)$ is surjective for $Z\in S$.

b) we have: $V_{Z_1\cap Z_2} = V_{Z_1} +  V_{Z_2}$ for
$Z_1,Z_2 \in S$.

c) the vector spaces $V_Z,\; Z\in S$  form a distributive
lattice of subspaces of $V$.}

(Recall that  a set of subspaces of vector space is called
a lattice if it is closed under the operations of sum
and intersections. A lattice is called distributive
if the distributiveness property : $V\cap (U+W)=V\cap U + V\cap W$
holds for any elements $V,U,W$ of the lattice.)

{\sf Proof:} a) We have an exact sequence: $H^0(Y,{\cal L})
\arrow  H^0({\cal L}\otimes {\cal O}_Z)\arrow H^1(Y, J^Y_Z\otimes
{\cal L})=0$.

b) An exact sequence of sheaves $0\arrow J_{Z_1 \cup Z_2}\arrow
J_{Z_1}\oplus J_{Z_2} \arrow  J_{Z_1\cap Z_2} \arrow 0$ remains exact
after tensoring with ${\cal L}$. The latter exact sequence
yields the cohomology sequence : $V_{Z_1}\oplus V_{Z_2}
\arrow V_{Z_1\cap Z_2} \arrow H^1(Y,J_{Z_1\cup Z_2}\otimes {\cal L})=0. $

c) It is obvious that $ V_{Z_1\cup Z_2}=V_{Z_1}\cap V_{Z_2}.$
Since any scheme in $S$ is reduced, the equality
$Z_1 \cup (Z_2\cap Z_3)=(Z_1\cup Z_2) \cap (Z_1 \cup Z_3)$
follows from the corresponding identity for sets. Thus c)
follows from b).\hfill\square

{\sf 2.2.} Recall that we consider varieties over an arbitrary
field $k$.

Fix a projective embedding  $\varphi : X\hookrightarrow {\Bbb P}^N$.
Denote $\FF :=\varphi ^*(\OO(1))$.
 Let $A:={\displaystyle \bigoplus _{n=0}^\infty  H^0(\FF ^{\otimes n})}$
be the projective coordinate ring of $X$, and $B_w$ be the projective
coordinate ring of a Schubert variety $X_w$.

 {\sf Theorem.} {\it a) $A$ is a Koszul quadratic algebra.

b) $B_w$ is a Koszul module over $A$.

c) $B_w$ is a Koszul quadratic algebra.}

(We refer to [BGSo] for the definition and properties
of Koszul modules over a Koszul algebra.)

{\sf Proof:}  We apply Lemma 2.1 in the following situation.
Take $Y$ to be $X^n$; the set  $S=S^n$ was desrcibed in 1.6. Put
$\LL := \FF ^{\boxtimes n}$. Let $i:X\hookrightarrow X^n$
 be the diagonal embedding and let $\Delta$ be its image.
By $\Delta _k$ we denote
the diagonal $x_k=x_{k+1}$. Note that $\Delta,\;\Delta _k$
 lie in $S^n$.

We have: $H^0(Y,\LL)=A_1 ^{\otimes n},\; H^0(X,i^*(\LL))=A_n$,
where $A_i$ is the $i$-th graded component of $A$. The restriction map
corresponds to multiplication in $A$.

Lemma 2.1  says that this map is surjective, i.e.
$A$ is generated by its first component.
By Lemma 2.1.ii) (applied succesively $(n-1)$ times) we see that the kernel
of multiplication map $A_1^{\otimes n} \arrow A_n$ is equal to
$\sum V_{\Delta _k}\subset H^0(X^n,\LL)=A_1 ^{\otimes n}$. It is obvious that
$V_{\Delta _k}=A_1^{\otimes k-1}\otimes I\otimes A_1^{\otimes n-k-1}$.
Thus we proved that $A$ is quadratic. Now the Koszul property follows
from Lemma 2.1.iii) by the
following well-known criterion:

{\sf 2.3.} {\it A quadratic algebra $A$ is Koszul iff  for all $n$
the lattice of subspaces
of $A_1^{\otimes n}$ generated by the spaces:
$A_1^{\otimes k-1}\otimes I\otimes A_1^{\otimes n-k-1},$
where $i=1,..,n-1$
 ("the lattice of relations") is distributive.}

Proof of b) is analogous. Add to the above the following notations:
$Z:=X_w\times X^{n-1},\;Z':=Z\cap \Delta$. Let $i_w,\;i'_w$
be the corresponding embeddings. We have: $H^0(i_w^*(\LL))=
(B_w)_n$.

Lemma 2.1.i) then implies that $B_w$ is generated
by $(B_w)_0$ as a module over $A$. Lemma 2.1 ii) implies that
the module of relations is generated by its component of degree 1.
The Koszul property of the module follows from Lemma 2.1.iii)
by the following criterion (see [BGSo]):

{\sf 2.4.} {\it Let $M$ be a graded module over a Koszul
quadratic algebra $A$. Assume that the multiplication map
$A\otimes M_0 \arrow M$ is surjective, and its kernel $I_M$
is generated by the first graded component $(I_M)_1$ (i.e.
$M$ is a quadratic module). Then $M$ is a Koszul module iff
  for all $n$ the lattice of subspaces of $M_0\otimes (A_1)^{\otimes n}$
generated by the spaces: $I_M\otimes  (A_1)^{\otimes n-1},$
$A_1^{\otimes k-1}\otimes I_M \otimes A_1^{\otimes n-k-1},$
where $i=1,..,n-1,$
is distributive.}

The statement c) follows from a), b) by the algebraic

{\sf 2.5. Lemma.} {\it Let $A$ be a be a Koszul algebra, and
$B$ be its quotient algebra. If $B$ is Koszul as a module over
$A$, then $A$ is a Koszul quadratic algebra.}

{\sf Proof:} Consider the spectral sequence:
$$E^2_{ij}=Tor^B_i(Tor ^A_j(B,k),k)\Rightarrow  Tor^A_{i+j}(k,k)$$

The action of $B$ on $Tor ^A_j(B,k)$ obviously factors through
$B/B^{>0}=k$. Thus  $E^2_{ij}=Tor ^A_j(B,k)\otimes
Tor ^B_i(k,k).$

The grading on $A$ and $B$ induces a grading on all the objects
under consideration.
Recall that an algebra $A$ is Koszul
 iff $Tor _i^A(k,k)$ is concentrated in  the graded
degree $i$. A module $M$ is Koszul iff
$Tor _i(M,k)$ is concentrated in  degree $i$.

Assume that $B$ is not Koszul. Then for some $i$ there exists
$s>i$ such that $Tor_i^B(k,k)$ has nonzero component of degree
$s$. Take the smallest $i$ with this property. Then we see that:
$E^2_{0,i}$ has nonzero component of degree $s>i$; and
$E^2_{p,q}$ is concentrated in degree $p+q$ for $p+q < i.$
(In the last assertion we used the Koszul property of $A$, and of $B$
as an $A$-module.) Since the differential $d_r :E^r _{ij}\arrow
E^r _{i+r,j-r-1}$ preserves the grading, we see that $d_r|_{(E^r_{i,0})_s}
=0$ for all $r.$ Thus $(E^\infty _{i,0})_s \not = 0$, and hence
 $Tor _i^A(k,k)$ has a nonzero part of degree $s$. But this contradicts
the Koszul property of $A.$\hfill \square

{\sf 2.5. Remark.} By analogy with the notion of Koszul module,
one can define a pair of graded modules $M,N$ over a Koszul algebra $A$
to be Koszul if $Tor _i^A(M,N)$ is concentrated in degree $i$ for
all $i>0.$ A criterion similar to 2.4 holds in this case. Our method
actually proves the following complement to Theorem 2.2:

{\it A pair (projective coordinate ring of a Schubert variety,
projective coordinate ring of an opposite Schubert variety)
is a Koszul pair of modules over $A.$}

{\sf 2.6. Remark.} One can deduce the exactness of Koszul
complex, or the homological criterion of Koszulity directly
from the cohomology vanishing 1.7 by using the sheaf-theoretic
interpretation of Koszul complex and bar-complex respectively.
(One of these possibilities was pointed out to me by M.Kapranov.)

{\bf 3. Proof of Proposition 1.6.} This section is an extension
of the arguments of [R] (cf. Theorem 3.5 there).

We fix some base field
$k$ of characteristic $p>0$.

{\sf 3.1. Lemma.} {\it
i) If $s$ is compatible with two subschemes $Y_1$, $Y_2$
then  it is compatible with there union and (scheme-theoretic)
intersection.

ii)Let $f: X\arrow X'$ be a morphism, such that $f_* (\OO _X)
= \OO _{X'}.$ Let $Y\subset X$ be a closed subscheme, $Y'$ be
its (scheme-theoretic) image. If $Y$ is compatibly split in
$X$, then  $Y'$ is compatibly split in
$X'.$

iii) Let $U\subset X$ is an open subscheme, which intersects
all components of a reduced closed subscheme $Y$. Assume that
a Frobenius splitting $s$ of $X$ is such that its restriction
to $U$ is compatible with $Y$. Then $s$ is compatible with $Y$.}

{\sf Proof:} see [R], 1.9; [MR], Proposition 4 or [R], 1.8;
[MR], Lemma 1 or [R] 1.7.\hfill \square

By Lemma 3.1.i) it sufficies to prove 1.6 for the case
$X=G/B$, where $B$ is a Borel subgroup. From  Lemma 3.1.ii) it follows that
it is enough to find a Frobenius splitting  $s\in Hom (Fr_*(\OO),
\OO)$ which is
compatible with each of the varieties a), b), c) of 1.6.
As was explained in [MR], we have:  $ Hom (Fr_*(\OO),
\OO)= H^n(Fr_*(\OO)\otimes K)^*=H^n ( Fr_*(Fr^*(K)))^*=
H^n ( Fr^*(K))^*= H^n (K^{\otimes p})^*=
Hom ( K^{\otimes p},K).$
(Here $K$ is the canonical class. The first and the last equality
are the Serre duality; the second is the projection formular, and the
third uses the fact that the finite morphism $Fr$ is acyclic.)
So we have to construct an appropriate element $s\in Hom ( K^{\otimes p},K)
= H^0(K^{\otimes 1-p}).$

Let $D\subset X^n$ be the union of all the varieties a), b), c) of 1.6.
Denote also by  $D^+\subset X$  the union of all Schubert varieties,
and by $D^-\subset X$  the union of all opposite Schubert varieties.

We have an isomorphism $\OO(D)\widetilde = K^{-1}.$
Indeed, since ${\displaystyle Pic(X^n)=\bigoplus _{i=1} ^n Pic(X)}$,
 it is enough to check that
for any $i$ the restrictions of  $\OO(D)$ and $K^{-1}$ to a fiber
 of the  projection
$p_i:X^n\arrow X^{n-1}$ (forgetting the $i$-th coordinate)
are isomorphic. So it reduces to a
well known isomorphism: $K_X^{-1}\widetilde = {\cal
O}(D^+ + D^-)$.

Hence there exists a section $\sigma\in H^0(X^n,K^{-1})$
whose divisor is equal to $D$. We put $s:=\sigma ^{p-1}.$

We claim that $s$ provides the desired Frobenius splitting.
In view of Lemma 3.1.iii) it is enough to check this locally.
Let us make this precise. Let $x\in X$ be the fixed point of
$B^+$, and $y\in X$ be the fixed point of $B^-$. By Lemma 3.1.iii),
the proof will be finished if we show the following:

{\sf 3.2. Lemma.}   {\it Let $s \in
  H^0(K^{-1})$ be a section defined on some neighbourhood
 of the point
$(x,..,x)\in X^n$ (respectively $(y,..,y)\in X^n$).
 Assume  that  the divisor of $s$ is equal to $D+D'$,
where $(x,..,x)\not \in {   supp}(D')$ (respectively,
 $(y,..,y)\not \in {   supp}(D')$). Then  $s^{p-1}$ gives a
Frobenius splitting of some neighbourhood of (x,...,x) which is compatible with
all varieties  a), b) (respectively, b), c)) of 1.6.}

To prove it we need the following key statement, which
is a consequence of the proof of Theorem 2 of [MR].

{\sf 3.3. Sublemma.}
 {\it Let  $\sigma \in
{   H}^0(K^{-1})$ be a section  defined on some neighbourhood of
 $x$ in $X$. Suppose that the divisor of $\sigma$ is equal
to $D^+ + D'$,
where $x\not \in {   supp}(D')$. Then $\sigma ^{p-1}$ gives a
Frobenius splitting of
some (possibly smaller) open neighbourhood of $x$ which is
compatible with all Schubert varieties.}\hfill \square

As an immediate corollary we get:

{\sf 3.4. Sublemma.}  {\it Let $\sigma \in
{   H}^0(K^{-1})$ be a section  defined on some neighbourhood
 of the point
$(x,..,x)\in X^n$. Suppose that  the divisor of
$\sigma$  is equal to $({\displaystyle \sum_{i}} X^i\times
D^+ \times X^{n-i})+D'$,
where $(x,..,x)\not \in {   supp}(D')$. Then  $\sigma^{p-1}$ provides a
Frobenius splitting of some neighbourhood of (x,...,x) which is compatible with
all  varieties $X^i\times X_w \times X^{n-i-1}$, where $X_w$
is a Schubert variety.}\hfill \square

Now we can deduce Lemma 3.2. First of all note, that there
exists an automorphism of $G$ taking $B^+$ to $B^-$, hence we can
restrict ourself to the case when $s$ is defined on a neighbourhood
of  $(x,..,x)$.

Let $N^-$ be the maximal nilpotent subgroup in $B^-$. We have an open embedding
$j:(N^-)^n\hookrightarrow
X^n$, such that: $j((g_1,...,g_n))=(g_1 x,...,g_n x)$, where for an element
$g\in N^-$ its action on $X$ is denoted by $g: y\arrow g y$.

 Define  an automorphism $a$ of $(N^-)^n$ by the formula
 $a:(g_1,...,g_n)\arrow (g_1, g_1^{-1}g_2,...,g_{n-1}^{-1}g_n)$. Then:
$jaj^{-1}((x,...,x))=(x,...,x)$. One can see
that $a$ sends $j^{-1}(D)$ to
 the divisor $j^{-1}({\displaystyle \sum_{i}} X^i\times
D^+\times X^{n-i})$. Hence Lemma 3.2 is equivalent to the
Sublemma 3.4. \hfill \square

{\bf Acknowledgements.} The statement a) of Theorem 2.2. was
communicated to me  as a conjecture by J. Donin. I would like
to express here my gratitude to him. I am glad for the
opportunity to thank L. Positsel'skii for various discussions
and explanations. Lemma 2.5 was pointed out to me by him.
Thanks for
encouragement  in writing the paper are due to V. Ginzburg and
J. Donin.

{\bf Remark.} During the final stage of preparation of the text
there appeared a paper by Mehta, Inamdar ``Frobenius splitting
of Schubert varieties and linear syzygies'' (American J. Math.,
116 (December, 1994)) on the same topic.
\hfill\break

\centerline{\bf References.}

 [BGS] A.Beilinson, V.Ginzburg, V.Schechtman {\it Koszul
duality}, JGP {\bf 5}, \# 3, (1988),  317-350.

[BGSo] A.Beilinson, V.Ginzburg, W.Soergel
{\it Koszul Duality Patterns in Representation Theory},
to appear (1995).

[BF] J.Backelin, R.Fr\"{o}berg.
{\it  Koszul algebras, Veronese subrings  and rings with linear
resolutions.}
 Rev. Roumaine Math. Pures Appl. {\bf 30}, \#2, (1985), p.85-97.

[EGA] A. Grothendieck, J.Dieudonn\' e {\it El\' ements de G\' eom\' etrie
Alg\' ebrique} III (partie 2), Publ. Math. IHES {\bf 17} (1963).

[LM] G. Lancaster,  J. Towber  {\it Representation-Functors
and Flag-algebras for the Classical Groups}, Jour. Alg. {\bf 59}
 (1979),  16-38.

[MR] V.B.Mehta, A.Ramanathan {\it Frobenius Splitting and Cohomology
Vanishing for Schubert Varieties}, Ann. Math. {\bf 122} (1985), 27-40.

 [PP] A.Polischuk, L.Positsel'skii {\it On  quadratic  algebras}, preprint,
Moscow, 1992.

[R] A. Ramanathan {\it Equations Defining Schubert Varieties
and Frobenius Splitting of Diagonals}, Publ. Math. IHES {\bf 65} (1987),
61-90.

\end{document}